\documentclass{jnmp01b}

\usepackage{amsmath}

\numberwithin{equation}{section}

\allowdisplaybreaks

\newcommand{\df}{\displaystyle}
\newcommand{\eps}{\epsilon}
\newcommand{\si}{\sigma}
\newcommand{\om}{\Omega}
\newcommand{\p}{\partial}
\newcommand{\la}{\langle}
\newcommand{\ra}{\rangle}

\newenvironment{namelist}[1]{%
\begin{list}{}
 {
   
   \settowidth{\labelwidth}{#1}
   \setlength{\leftmargin}{1.1\labelwidth}
  }
 }{%
\end{list}}
\newcommand{\bi}{\begin{namelist}}
\newcommand{\ei}{\end{namelist}}

\newtheorem{prop}{Proposition}[section]
\newtheorem{theorem}{Theorem}[section]
\newtheorem{lemma}{Lemma}[section]
\theoremstyle{definition}
\newtheorem{definition}{Definition}[section]
\newtheorem{remark}{Remark}

\makeatletter
\DeclareRobustCommand{\primfrac}[1]{%
  \PackageWarning{amsmath}{%
Foreign command \@backslashchar#1; %
\protect\frac\space or \protect\genfrac\space should be used instead%
  }
  \global\@xp\let\csname#1\@xp\endcsname\csname @@#1\endcsname
  \csname#1\endcsname
}
\makeatother

\begin{document}

\renewcommand{\evenhead}{N.\ Svanstedt}
\renewcommand{\oddhead}{Correctors for the Homogenization of 
Monotone Parabolic Operators}


\thispagestyle{empty}

\begin{flushleft}
\footnotesize \sf
Journal of Nonlinear Mathematical Physics \qquad 2000, V.7, N~3,
\pageref{firstpage}--\pageref{lastpage}.
\hfill {\sc Letter}
\end{flushleft}

\vspace{-5mm}

\copyrightnote{2000}{N.\ Svanstedt}

\Name{Correctors for the Homogenization of Monotone Parabolic Operators}

\label{firstpage}

\Author{Nils SVANSTEDT}

\Adress{Department of Mathematics \\
        Chalmers University and G\"{o}teborg University \\
        S-412 96 G\"{o}teborg, SWEDEN}

\Date{Received October 11, 1999; Revised March 17, 2000; 
Accepted May 5, 2000}

\begin{abstract}
\noindent
In the homogenization of monotone
parabolic partial differential equations
with oscillations in both the space and time variables
the gradients converges only weakly in $L^p$. In the
present paper we
construct a family of correctors, such that, up to a remainder which
converges to zero strongly in $L^p$, we obtain strong convergence of   
the gradients in $L^p$.
\end{abstract}


\section{Introduction}

In [8] the asymptotic behaviour
(as $\eps\to{0}$) of the solutions $u_{\eps}$ to a sequence of
initial-boundary value problems of the form \\
\begin{equation}
\left\{ \begin{array}{l}
{\df {\p u_{\eps}\over \p t} - {\rm div}(a({x\over \eps},{t\over \eps^{\mu}},Du_{\eps})) = f}\;
\mbox{\rm  in }\om\times]0,T[, \\[2ex]
u_{\eps}(x,0) = u_0(x), \\[1ex]
u_{\eps}(x,t) = 0\;\mbox{\rm  in }\p\om\times]0,T[, 
\end{array} \right.
\end{equation}
is studied. Here $\om$ is an open bounded set in ${\bf R}^{\rm N}$,
T is a positive real number, 
$2\leq{p}<\infty$ and $\mu>0$. Under the assumption that
$a({x\over \eps},{t\over \eps^{\mu}},Du_{\eps})$ is $\eps-$ 
and $\eps^\mu$-periodic in the first and second variable, respectively,
it is proved that
\begin{gather*}
u_{\eps}\to{u}\;\mbox{\rm weakly in }
\;{\rm L}^p(0,{\rm T};{\rm W}^{1,p}_0(\om)), \\
\textstyle
a({x\over \eps},{t\over \eps^{\mu}},Du_{\eps})\to{b(Du)}
\;\mbox{\rm weakly in }
\;{\rm L}^{p'}(0,{\rm T};{\rm L}^{p'}(\om;{\bf R}^{\rm N})),
\end{gather*}
where $1/p + 1/p' = 1$ and
where $u$ denotes the unique solution to
\begin{equation}
\left\{ \begin{array}{l}
{\df {\p u\over \p t} - {\rm div}(b(Du)) = f}\;\;{\rm  in }\;\;\om\times]0,T[, 
\\[2ex]
u(x,0) = u_0(x), \\[1ex]
u(x,t) = 0\;\mbox{\rm  in }\p\om\times]0,T[, 
\end{array} \right.
\end{equation}
where the limit map $b$ in (1.2) only depends on the sequence
$(a(\cdot/\eps,\cdot/\eps^\mu,\xi))$
and on $\mu$ and where $b$ is different for $0<{\mu}<{2}$, $\mu=2$ and $\mu>{2}$, respectively.
This result implies that
\[
Du_{\eps}(x,t) = Du(x,t) + r_{\eps}(x,t),
\]
where the remainder $r_{\eps}$ converges to zero only weakly in
${\rm L}^{p}(0,{\rm T};{\rm L}^{p}(\om;{\bf R}^{\rm N}))$. The purpose of the present
paper is to construct a family of correctors $(p_\eps)=(p_{\eps}(x,t,\xi))$ such that
\[
p_{\eps}(\cdot,\cdot,\xi)\to\xi\mbox{ weakly in }\;
{\rm L}^{p}(0,{\rm T};{\rm L}^{p}(\om;{\bf R}^{\rm N}))
\]
for every $\xi\in{\bf R}^{\bf N}$ and
\[
Du_{\eps}(x,t) = p_{\eps}(x,t,(M_{\eps}Du)(x,t)) + r_{\eps}(x,t),
\]
where the remainder $r_{\eps}$ converges to zero strongly in
${\rm L}^{p}(0,{\rm T};{\rm L}^{p}(\om;{\bf R}^{\rm N}))$ and where
$(M_{\eps})$ is a sequence of linear operators
on ${\rm L}^{p}(0,{\rm T};{\rm L}^{p}(\om;{\bf R}^{\rm N}))$ which
converges strongly to the identity map on
${\rm L}^{p}(0,{\rm T};{\rm L}^{p}(\om;{\bf R}^{\rm N}))$.
The results presented in this article are rather technical and involves numerous
estimates on small $\eps$-cubes. But the implications from Theorem 2.1 are important.
In particular for computational modeling of (1.1)
and (1.2) since it implies strong convergence
of the gradients in the energy norm. In a simplified way we can say that
the improvement of the convergence lies in the fact that the
local behaviour on the $\eps$-cubes
are added to the homogenized solution. Heuristically this amounts to adding the
second term in an asymptotic expansion, see [1].
The corrector problem was first studied in [1] for linear
elliptic and parabolic problems. For a careful study of
linear parabolic problems we refer to
[3]. See also [2] and [6].
The extension to the
monotone elliptic case is performed in [4]. The present work is very much inspired
by the methods developed in [4].
The paper is organized as follows: Section 2 contains
some preliminaries and in Section 3 we present the main theorem (Theorem 3.1).
In Section 4 we collect some useful estimates for the correctors and in 
Section 5 we give the proof of Theorem 3.1.

\section{Preliminaries}

Throughout this paper we will denote by $\om$
a bounded open set in ${\bf R}^{\rm N}$ and
we will let ${\rm V} = {\rm W}^{1,p}_0(\om)$,
with norm $\|u\|_V^p = \int_{\om}|Du|^p\,dx$ and
${\rm V'} = {\rm W}^{-1,p'}(\om)$.

We consider the
evolution triple
\[
{\rm V}\subseteq{\rm L}^2(\om)\subseteq{\rm V'},
\]
with dense embeddings. Further, for positive
real-valued T and for $2\leq{p}<\infty$, we define
${\cal V} = {\rm L}^p(0,{\rm T};{\rm V})$ and
${\cal V'} = {\rm L}^{p'}(0,{\rm T};{\rm V'})$, where $1/p + 1/p' = 1$ and
the corresponding evolution triple
\[
{\cal V}\subseteq{{\rm L}^2(]0,{\rm T}[\times{\om})}\subseteq{\cal V'}
\]
also with dense embeddings where the duality pairing
$\la\cdot,\cdot\ra_{\cal V}$
between ${\cal V}$ and ${\cal V}'$ is
given by
\[
\la{f},u\ra_{\cal V} =
\int^T_0\la{f}(t),u(t)\ra_V\,dt,\;{\rm for}\;u\in{\cal V},\;
f\in{\cal V}'.
\]

Given $u_0\in{\rm L}^2(\om)$, the space ${\cal W}_0$ is defined as
\[
{\cal W}_0 = \{v\in{\cal V}:v'\in{\cal V'}\;{\rm and }\;v(0)
=u_0\in{\rm L}^2(\om)\}.
\]
Here $v'$ denotes the time derivative of $v$, which
is to be taken in distributional sense. 
Moreover, we define
\[
{\cal U} = {\rm L}^p(0,{\rm T};{\rm L}^p(\om;{\bf R}^{\rm N}))
\mbox{ and  }\;
{\cal U}' = {\rm L}^{p'}(0,{\rm T};{\rm L}^{p'}(\om;{\bf R}^{\rm N})).
\] 
with the duality pairing
\[
\la{u},v\ra_{\cal U} = \int^T_0\int_{\om}(u,v)\,dxdt,\;\mbox{ for }u\in{\cal U}'
\mbox{ and }v\in{\cal U},
\] 
where $(\cdot,\cdot)$ denotes the scalar product in ${\bf R}^{\rm N}$. By $|\cdot|$
we understand the usual Euclidean norm in ${\bf R}^{\rm N}$
and by $m(\cdot)$ we understand the Lebesgue measure. Moreover, by
$(\eps)$ we understand a sequence of positive
real numbers tending to $0^{+}$.\\

Let $Y=]0,1[^{\rm N}$ be the unit cube in ${\bf R}^{\rm N}$ and let
$Y\times{T_0}=]0,1[^{\rm N}\times{]0,1[}$ be the unit cube in
${\bf R}^{\rm N}\times{\bf R}_{+}$.
\begin{definition} 
We say that a function $u:{\bf R}^{\rm N}\times{]0,{\rm T}[}\to
{\bf R}$ is $Y$-{\em periodic} if $u(x+e_i,t) = u(x,t)$ for every
$x\in{\bf R}^{\rm N}$, $t\in{]0,{\rm T}[}$ and for every
$i=1,\ldots,\,{\rm N}$, where $(e_i)$ is the canonical basis of
${\bf R}^{\rm N}$. Further, we say that a function $u:{\bf R}^{\rm N}\times{\bf R}_{+}\to
{\bf R}$ is $Y\times{T_0}$-{\em periodic} if $u(x+e_i,t) = u(x,t)
= u(x,t +1)$ for every
$x\in{\bf R}^{\rm N}$, $t\in{\bf R}_{+}$ and for every
$i=1,\ldots,\,{\rm N}$.
\end{definition}

We consider the following
spaces of periodic functions:
\[
{\rm V}_{\sharp,Y} = \{u\in{\rm W}^{1,p}_{\rm loc}({\bf R}^{\rm N}):u
\mbox{ is Y-periodic and has mean value zero over Y}\},
\]
and
\[
{\cal V}_{\sharp,Y\times\ {T_0}}
= \{u\in{\rm L}^{p}_{\rm loc}({\bf R}_+;{\rm V}_{\sharp,Y}):u
\mbox{ is  }T_0-periodic\}.
\]
\begin{definition} 
Given $0<\alpha\leq{1}$, $2\leq{p}<\infty$ and three positive real
constants $c_0$, $c_1$ and $c_2$ we define the class 
$S_{\sharp,Y\times{T_0}} = S_{\sharp,Y\times{T_0}}(c_0,c_1,c_2,\alpha)$
of maps
\[
a:{\bf R}^{\rm N}\times{\bf R}_{+}\times{\bf R}^{\rm N}\to{\bf R}^{\rm N},
\]
such that
\bi{(iii)xx}
\item[{\rm (i)}] $a(\cdot,\cdot,\xi)$ is $Y\times{\tau_0}$-periodic
for every $\xi\in{\bf R}^{\rm N}$,
\item[{\rm (ii)}] $|a(y,\tau,0)|\leq{c_0}$ a.e in ${\bf R}^{\rm N}\times{\bf R}_{+}$,
\item[{\rm (iii)}] $a(\cdot,\cdot,\xi)$ is Lebesgue measurable
for every $\xi\in{\bf R}^{\rm N}$,
\item[{\rm (vi)}] $|(a(y,\tau,\xi_1) - a(y,\tau,\xi_2)|
\leq{c_1}(1 + |\xi_1| + |\xi_2|)^{p-1-\alpha}|\xi_1 - \xi_2|^{\alpha}$, a.e. in
${\bf R}^{\rm N}\times{\bf R}_{+}$
for all $\xi_1, \xi_2\in{\bf R}^{\rm N}$,
\item[{\rm (v)}] $(a(y,\tau,\xi_1) - a(y,\tau,\xi_2),\xi_1 - \xi_2) \geq{c_2}|\xi_1 - \xi_2|^p$,
a.e. in ${\bf R}^{\rm N}\times{\bf R}_{+}$
for all $\xi_1, \xi_2\in{\bf R}^{\rm N},\;\xi_1\neq{\xi_2}$.
\ei
\end{definition}

We recall some results for maps $a\in{S_{\sharp,Y\times{T_0}}}$:
\begin{prop}
Suppose that $a\in{S_{\sharp,Y\times{T_0}}}$. Then,
for every $f\in{\cal V'}$ and for every $\eps>0$,
(1.1) possesses a unique solution $u_{\eps}\in{\cal W}_0\cap L^\infty(0,T;L^2(\om))$.
\end{prop}
\begin{proof}
See e.g. [10].
\renewcommand{\qed}{}
\end{proof}
\begin{prop} Let us put
$a_\eps(\cdot,\cdot,\xi) = a({\cdot\over \eps},{\cdot\over \eps^\mu},\xi)$.
Suppose that $a\in{S}_{\sharp,Y\times{T_0}}$. Then, for every
$f\in{\cal V'}$, the solutions $u_{\eps}$ to (1.1) satisfy
\[
\left. \begin{array}{l}
u_{\eps}\to{u}\;\mbox{\rm weakly in }
\;{\cal W}_0, \\\\
a_{\eps}(x,t,Du_{\eps})\to{b(Du)}
\;\mbox{\rm weakly in }
\;{\cal U'},
\end{array} \right.
\]
where $u$ is the unique solution to the following parabolic problem:
\begin{equation}
\left\{ \begin{array}{l}
u' - {\rm div}(b(Du)) = f\mbox{ in }\om\times]0,T[ \\
u\in{\cal W}_0. 
\end{array} \right.
\end{equation}
Moreover, for a 
fixed vector $\xi\in \bf R^{\rm N}$:\\
\begin{equation}
b(\xi) = \int_{\tau_0}\int_Ya(y,\tau,Dv(y,\tau)+\xi)\,dyd{\tau},
\end{equation}
where $v$ depends on $\xi$ and $\mu$. For
$0<{\mu}<{2}$, $v=v(y,\tau)$ is the unique solution
to the parameter-dependent elliptic
problem: 
\begin{equation}
\left\{ \begin{array}{l}
- {\rm div}(a(y,\tau,Dv(y,\tau)+\xi)) = 0,\\
v(\cdot,\tau)\in{\rm V}_{\sharp,Y},\,\tau\geq 0. 
\end{array} \right.
\end{equation}
For $\mu=2$, $v=v(y,\tau)$ is the unique solution, to the parabolic
problem: 
\begin{equation}
\left\{ \begin{array}{l}
v'- {\rm div}(a(y,\tau,Dv(y,\tau)+\xi)) = 0,\\
v\in{\cal V}_{\sharp,Y\times{T_0}}. 
\end{array} \right.
\end{equation}
For $\mu>2$, finally $v=v(y)$ is the unique solution to the elliptic
problem: 
\begin{equation}
\left\{ \begin{array}{l}
- {\rm div}({\tilde {a}}(y,Dv(y)+\xi)) = 0,\\
v\in{\rm V}_{\sharp,Y}, 
\end{array} \right.
\end{equation}
where
\begin{equation}
{\tilde {a}}(y,\xi) = \int_{\tau_0}a(y,\tau,\xi)d\tau.
\end{equation}
\end{prop}
\begin{proof}
We refer to [8].
\renewcommand{\qed}{}
\end{proof}
\begin{remark} 
By the estimates (4.19) and (4.21)
in the proof of Theorem 3.1 in [9]
it follows that the homogenized map $b$ satisfies the estimates
\begin{gather*}
|b(\xi_1)-b(\xi_2)|\leq{C}(1+|\xi_1|+|\xi_2|)^{p-1-\gamma}|\xi_1-\xi_2|^{\gamma}
\\[1ex]
(b(\xi_1)-b(\xi_2),\xi_1-\xi_2)\geq{c_2}|\xi_1-\xi_2|^p
\end{gather*}
for every $\xi_1,\,\xi_2\in{\bf R}^{\rm N}$ where $\gamma =
\alpha/(p-\alpha)$.
\end{remark}

We close this section by stating some different Meyers type estimates which will be needed in the
proof of the main corrector result, Theorem 3.1.
\begin{prop} 
Suppose that $a\in{S_{\sharp,Y\times{T_0}}}$.
Let $u$ be the solution to the problem
\[
\left\{ \begin{array}{l}
- {\rm div}({\tilde {a}}(x,Du) = 0,\\
u\in W^{1,p}(\om), 
\end{array} \right.
\] 
Then there exists a constant $\eta>0$ such that
$u\in W^{1,p+\eta}({\tilde \om})$
for every open set ${\tilde \om}\subset\subset\om$. Moreover
\[
\|u\|_{W^{1,p+\eta}({\tilde \om})}\leq C\|u\|_{W^{1,p}({\om})}.
\]
\end{prop}
\begin{proof}
We refer to Theorem 1 in [7].
\renewcommand{\qed}{}
\end{proof}
\begin{remark}
Considered as a function constant in $t\geq{0}$, the
function $u$ above also satisfies the estimate
\[
\|u\|_{L^{p+\eta}(0,T;W^{1,p+\eta}({\tilde \om}))}\leq C\|u\|_{L^p(0,T;W^{1,p}({\om}))}.
\]
\end{remark}
\begin{prop}
Suppose that $a\in{S_{\sharp,Y\times{T_0}}}$
and in addition satisfies
\begin{equation}
|a(x,t,\xi)-a(x,s,\xi)|\leq\omega(t-s)(1+|\xi|^{p-1})
\end{equation}
for all
$t,s\in ]0,T[$, all $\xi\in{\bf R}^{\rm N}$ and a.e. $x\in \om$, where
$\omega$ is the modulus of continuity.
Let $u(\cdot,t)$, $t\in ]0,T[$, be the solution to the
parameter dependent elliptic problem
\[
\left\{ \begin{array}{l}
- {\rm div}(a(x,t,Du)) = 0,\\
u(\cdot,t)\in W^{1,p}(\om), 
\end{array} \right.
\] 
Then there exists a constant $\eta>0$ such that, for every $t\in ]0,T[$,
$u(\cdot,t)\in W^{1,p+\eta}({\tilde \om})$
for every open set ${\tilde \om}\subset\subset\om$. Moreover
\begin{equation}
\|u(\cdot,t)\|_{W^{1,p+\eta}({\tilde \om})}\leq C\|u(\cdot,t)\|_{W^{1,p}({\om})}.
\end{equation}
Further, let $\delta\subset\subset]0,T[$. 
The gradient $Du$ of the function $u$ above also satisfies the estimate
\begin{equation}
\|Du\|_{L^{p+\eta}(\delta;L^{p+\eta}({\tilde \om};{\bf R}^{\rm N}))}\leq
{C}\|Du\|_{L^p(0,T;L^p({\om};{\bf R}^{\rm N}))}.
\end{equation}
\end{prop}
\begin{proof}
The estimate (2.8) is a consequence of Theorem 1 in [7], if we
take (2.7) into account. By using the coercivity of $a$ and the H\"{o}lder
inequality we get
\begin{gather*}
\int_\delta\|Du(\cdot,t+\eps)-Du(\cdot,t)\|^p_{L^p(\om)}\,dt  \\ 
\quad\leq c_2\int_\delta\int_\om(a(x,t+\eps,Du(x,t+\eps))-a(x,t,Du(x,t)),Du(x,t+\eps) - Du(x,t))\,dxdt  \\ 
\quad\leq c_2(\int_\delta\|a(\cdot,t+\eps,Du(\cdot,t+\eps))-a(\cdot,t,Du(\cdot,t))
\|^{p'}_{L^{p'}(\om)}\,dt)^{1/p'}  \\ 
\quad\qquad\times (\int_\delta\|Du(\cdot,t+\eps)-Du(\cdot,t)\|^p_{L^p(\om)}\,dt)^{1/p}. 
\end{gather*}
By (2.7) we obtain, using the Minkowski inequality and the boundedness of $\om$,
\begin{gather*}
(\int_\delta\|Du(\cdot,t+\eps)-Du(\cdot,t)\|^{p}_{L^p(\om)}\,dt)^{p'} \\
\qquad\leq c_2(\int_\delta\|\omega(\eps)(1+|\max\{Du(\cdot,t+\eps),Du(\cdot,t)\}|^{p-1}
\|^{p'}_{L^{p'}(\om)}\,dt)^{1/p'} \\
\qquad\leq c_2(\int_\delta|\omega(\eps)|((m(\om))^{p'}+
\|\max\{Du(\cdot,t+\eps),Du(\cdot,t)\}\|^{p}_{L^{p}(\om)})\,dt)^{1/p'}.
\end{gather*}
For $\eps$ small enough $\omega(\eps)\leq 1$ and (2.8) implies that
\[
\sup_{t\in\delta}\|\max\{Du(\cdot,t+\eps),Du(\cdot,t)\}\|^{p}_{L^{p}(\om)}\leq C,
\]
where $C$ is independent of $\eps$. Therefore
\[
\int_\delta\|Du(\cdot,t+\eps)-Du(\cdot,t)\|^{p}_{L^p(\om)}\,dt
\leq C\int_\delta|\omega(\eps)|\,dt,
\]
which tends to zero
as $\eps\to 0$, by the dominated convergence theorem.
The estimate (2.9) now readily follows by the continuity of $Du$ with respect to $t$.
\end{proof}
\begin{prop}
Suppose that $a\in{S_{\sharp,Y\times{T_0}}}$.
Let $u\in L^p(0,T;W^{1,p}(\om))\cap L^\infty(0,T;L^2(\om))$ be the solution to the problem
\[
u'- {\rm div}(a(x,t,Du) = 0.
\]
Let ${\tilde \om}$ be defined as above and let $\delta\subset\subset]0,T[$. 
The gradient $Du$ of the function $u$ above also satisfies the estimate
\[
\|Du\|_{L^q(\delta;L^q({\tilde \om};{\bf R}^{\rm N}))}\leq
{C}\|Du\|_{L^p(0,T;L^p({\om};{\bf R}^{\rm N}))}.
\]
for any $q\in]1,\infty[$.
\end{prop} 
\begin{proof}
We refer to Lemma 2.2 and Remark 7.4 of [5].
\renewcommand{\qed}{}
\end{proof}

\section{The main result}

In this section we state the main corrector result which we indicated in the
previous sections. We start out by defining a sequence
$(M_{\eps})$ of approximations of the identity map on
${\cal U}$. For
$i\in{\bf Z}^{\rm N}$ and  $j\in{\bf Z}$ we consider the translated
images $Y^i_{\eps} = \eps(i+Y)$ and ${T^j_{0,\eps}} =
\eps^\mu(j + T_0)$. Take
$\varphi\in{\cal U}$.
We define the function
\[
M_{\eps}\varphi:{\bf R}^{\rm N}\times{\bf R}\to{\bf R}^{\rm N}
\]
by
\begin{equation}
(M_{\eps}\varphi)(x,t) = \sum_{i\in{I_{\eps}}}\sum_{j\in{J_{\eps}}}
\chi_{Y^i_{\eps}}(x)\chi_{T^j_{0,\eps}}(t){1\over m(Y^i_{\eps}\times
{T^j_{0,\eps}})}\int_{T^j_{0,\eps}}\int_{Y^i_{\eps}}\varphi(y,\tau)
\,dyd{\tau},
\end{equation}
where
\[
I_{\eps} = \{i\in{\bf Z}^{\bf N}:Y^i_{\eps}\in{\om}\}\mbox{ and }
J_{\eps} = \{j\in{\bf Z}:{T}^j_{0,\eps}\in{]0,{\rm T}[}\}
\]
and $\chi_A$ denotes the characteristic function of the measurable set $A$.
It is well-known that
\begin{equation}
M_{\eps}\varphi\to\varphi\;\mbox{ strongly in }\;
{\cal U}.
\end{equation}
By the Jensen's inequality we also have
\begin{equation}
\|M_{\eps}\varphi\|_{\cal U}\leq\|\varphi\|_{\cal U}
\end{equation}
for all $\varphi\in{\cal U}$.

Let us also define the $Y\times{T_0}$-periodic function
\[
p:{\bf R}^{\rm N}\times{\bf R}\times{\bf R}^{\rm N}\to{\bf R}^{\rm N},
\]
which depends on $\mu$, by
\begin{equation}
p(x,t,\xi) = \xi + Dv(x,t),
\end{equation}
where $v$ is the solution to the auxiliary local problem (2.3),
(2.4) or (2.5) for $0<\mu<2$, $\mu=2$ and $\mu>2$,
respectively. It follows that
the function
\[
p_{\eps}:{\bf R}^{\rm N}\times{\bf R}\times{\bf R}^{\rm N}\to{\bf R}^{\rm N}
\]
defined by
\begin{equation}
p_{\eps}(x,t,\xi) = \xi + Dv({x\over \eps},{t\over \eps^\mu}),
\end{equation}
is $\eps{Y}$-periodic in $x$ and $\eps^\mu{T_0}$-periodic in $t$. This means that
\[
\int_{T_0}\int_{Y}p(x,t,\xi)\,dxdt = \xi
\]
and that
\begin{equation}
p_{\eps}(\cdot,\cdot,\xi)\to\xi\;\mbox{ weakly in }\;{\cal U}.
\end{equation}
Thus, the homogenized map $b$ can be expressed as
\begin{equation}
b(\xi) = \int_{T_0}\int_{Y}a(x,t,p(x,t,\xi))\,dxdt.
\end{equation}
Moreover, we have
\begin{gather}
\int_{T_0}\int_{Y}(a(x,t,p(x,t,\xi)),p(x,t,\xi))\,dxdt \notag\\
\qquad =\int_{T_0}\int_{Y}(a(x,t,p(x,t,\xi)),\xi)\,dxdt = (b(\xi),\xi).
\end{gather}

The following correctors result is the main result of this paper:
\begin{theorem}
Suppose that $a\in{S_{\sharp,Y\times{T_0}}}$. For the case
$0<\mu< 2$ we also suppose that $a$ satisfies (2.7).
Moreover, suppose that $f\in{\cal V}'$ and let $u_{\eps}$ be the
solutions to {\rm (1.1)} and let $u$ be the solution to {\rm (1.2)}.
Then, we have
\begin{equation}
Du_{\eps} = p_{\eps}(\cdot,\cdot,M_{\eps}Du) + r_{\eps},
\end{equation}
where $p_{\eps}$ is defined by {\rm (3.5)} and where
\[
r_{\eps}\to{0},\;\mbox{ strongly in }\;{\cal U}.
\]
\end{theorem}
\begin{remark} 
Recall that $p_\eps$ is entirely different for the three
cases $0<\mu<2$, $\mu=2$ and $\mu>2$, respectively.
\end{remark}

\section{Some estimates for the family of correctors}

In this section we present some estimates
for the family $(p_\eps)$ of correctors.
To a large extent the proofs will follow by minor
modifications of the proofs of similar lemmas by Dal Maso and
Defranceschi in [4]. Therefore we refer to their paper for
complete details and present here only proofs of parts which require
more modifications.
\begin{lemma} For any vector $\xi\in{\bf R}^{\rm N}$ we have
\begin{equation}
\|p_{\eps}(\cdot,\cdot,\xi)\|^p_{L^p(T_{0,\eps};L^p(Y_{\eps};{\bf R}^{\rm N}))}
\leq{C}(1+|\xi|^p)m(Y_{\eps}\times{T_{0,\eps}}),
\end{equation}
where the constant $C$ depends only on ${\rm N}$, $p$, $c_0$, $c_1$
and $c_2$.
\end{lemma}
\begin{lemma} There exist $\eta>{0}$ and $C>{0}$, which depends only on
${\rm N}$, $p$, $c_0$, $c_1$
and $c_2$, such that
\begin{equation}
\|p_{\eps}(\cdot,\cdot,\xi)\|^{p+\eta}_{L^{p+\eta}
(T_{0,\eps};L^{p+\eta}(Y_{\eps};{\bf R}^{\rm N}))}
\leq{C}(1+|\xi|^{p+\eta})m(Y_{\eps}\times{T_{0,\eps}}),
\end{equation}
for every $\xi\in{\bf R}^{\rm N}$.
\end{lemma}
\begin{proof}
By referring to the Meyers estimates in Propositions 2.3, 2.4 and 2.5
the proof is analogous as the proof of Corollary 3.3 in [4].
\end{proof}
\begin{lemma} For every $\xi_1,\xi_2$ in ${\bf R}^{\rm N}$
we have
\begin{gather}
\|p_{\eps}(\cdot,\cdot,\xi_1)-p_{\eps}(\cdot,\cdot,\xi_2)\|
^p_{L^p(T_{0,\eps};L^p(Y_{\eps};{\bf R}^{\rm N}))} \notag\\[1ex]
\qquad\leq {C}(1+|\xi_1|^p+|\xi_2|^p)^{(p-1-\alpha)/(p-\alpha)}
|\xi_1-\xi_2|^{p/(p-\alpha)}m(Y_{\eps}\times{T_{0,\eps}}), 
\end{gather}
where the constant $C$ depends only on ${\rm N}$, $p$, $\alpha$, $c_0$, $c_1$
and $c_2$.
\end{lemma}
\begin{lemma}
Let $\varphi\in{\cal U}$ and consider a simple function $\Psi$ given by
\[
\Psi(x,t) = \sum^{m}_{k=1}c_k\chi_{\om_k}(x)\chi_{\delta_k}(t),
\]
with $c_k\in{\bf R}^{\rm N}\setminus\{0\}$, $\om_k\subset\subset\om$,
$\delta_k\subset\subset{]0,{\rm T}[}$, $m(\p\om_k) = m(\p\delta_k) = 0$
and $(\om_k\cap\om_l)\times(\delta_k\cap\delta_l)=\phi$ for $k\neq{l}$.
Then
\begin{gather}
\limsup_{\eps\to{0}}\|p_{\eps}(\cdot,\cdot,M_{\eps}\varphi)-
p_{\eps}(\cdot,\cdot,\Psi)\|_{\cal U} \notag \\
\qquad\leq C(m(\om\times{]0,{\rm T}[})+\|\varphi\|_{\cal U}+\|\Psi\|_{\cal U})^
{(p-1-\alpha)/(p-\alpha)}\|\varphi-\Psi\|^{1/(p-\alpha)}_{\cal U},
\end{gather}
where $C$ depends only on ${\rm N}$, $p$, $\alpha$, $c_0$, $c_1$
and $c_2$.
\end{lemma}
\begin{proof}
Put $\om_0=\om\setminus\cup^m_{k=1}\om_k$,
$\delta_0={]0,{\rm T}[}\setminus\cup^m_{k=1}\delta_k$ and $c_0=0$. Then
we have
\[
\Psi(x,t) = \sum^{m}_{k=0}c_k\chi_{\om_k}(x)\chi_{\delta_k}(t).
\]

For every $\eps>{0}$ we denote by $\om_{\eps}\times\delta_{\eps}$
the union of all closed cubes
${\overline Y}^i_{\eps}\times{\overline T}^j_{0,\eps}$ such that
$Y^i_{\eps}\subset\om$ and $T^j_{0,\eps}\subset{]0,{\rm T}[}$.
For $k=0,\,1,\ldots,\,m$ we define the sets
\[
I^k_{\eps}=\{i\in{I_{\eps}}:Y^i_{\eps}\subset\om_k\},\;\;
J^k_{\eps}=\{j\in{J_{\eps}}:T^i_{0,\eps}\subset\delta_k\},
\]
and
\begin{gather*}
{\tilde {I}}^k_{\eps}=\{i\in{I_{\eps}}:
Y^i_{\eps}\cap\om_k\neq\phi,\;Y^i_{\eps}\setminus\om_k\neq\phi\},\\
{\tilde {J}}^k_{\eps}=\{j\in{J_{\eps}}:
T^i_{0,\eps}\cap\delta_k\neq\phi,\;T^i_{0,\eps}\setminus\delta_k\neq\phi\}.
\end{gather*}

Further, we define $E^{i,j,k}_{\eps}$ as the union of all closed cubes
${\overline Y}^i_{\eps}\times{\overline T}^j_{0,\eps}$ with
$i\in{I^k_{\eps}}$ and $j\in{J^k_{\eps}}$, and we define
${\tilde {E}}^{i,j,k}_{\eps}$ as the union of all closed cubes
${\overline Y}^i_{\eps}\times{\overline T}^j_{0,\eps}$ with
$i\in{\tilde {I}}^k_{\eps}$ and $j\in{\tilde {J}}^k_{\eps}$.
If we choose $\eps$ small enough, then, for $k\neq{0}$,
$\om_k\times\delta_k\subseteq\om_{\eps}\times\delta_{\eps}$ according
to (3.1). Thus, the definition of $\Psi$ yields
\begin{gather*}
\|p_{\eps}(\cdot,\cdot,M_{\eps}\varphi)
-p_{\eps}(\cdot,\cdot,\Psi)\|^p_{\cal U}
=\int_{\delta_{\eps}}\int_{\om_{\eps}}
|p_{\eps}(x,t,M_{\eps}\varphi)-p_{\eps}(x,t,\Psi)|^p\,dxdt  \\
\qquad\leq \sum_{k=0}^m\int_{E^{i,j,k}_{\eps}}
|p_{\eps}(x,t,M_{\eps}\varphi)-p_{\eps}(x,t,c_k)|^p\,dxdt \\*
\qquad\quad{}+\sum_{k=0}^m\int_{{\tilde {E}}^{i,j,k}_{\eps}}
|p_{\eps}(x,t,M_{\eps}\varphi)-p_{\eps}(x,t,c_k)|^p\,dxdt. 
\end{gather*}
Let us put
\[
\theta^{i,j}_{\eps} = {1\over m(Y^i_{\eps}\times{T}^j_{0,\eps})}
\int_{T^j_{0,\eps}}\int_{Y^i_{\eps}}\varphi(x,t)\,dxdt.
\]
A repeated application of the H\"{o}lder's and the Jensen's inequalities
yields, according to Lemma 4.3,
\begin{gather*}
\|p_{\eps}(\cdot,\cdot,M_{\eps}\varphi)
-p_{\eps}(\cdot,\cdot,\Psi)\|^p_{\cal U} \\
\quad\leq{C}\sum_{k=0}^m
\left(\sum_{j\in{J}^k_{\eps}}\sum_{i\in{I}^k_{\eps}}
(1+|\theta^{i,j}_{\eps}|^p+|c_k|^p)^{(p-1-\alpha)/(p-\alpha)}
|\theta^{i,j}_{\eps}-c_k|^{p/(p-\alpha)}
m(Y^i_{\eps}\times{T^j_{0,\eps}})\right) \\
\quad{}+{C}\sum_{k=0}^m
\left(\sum_{j\in{\tilde{J}}^k_{\eps}}\sum_{i\in{\tilde {I}}^k_{\eps}}
(1+|\theta^{i,j}_{\eps}|^p+|c_k|^p)^{(p-1-\alpha)/(p-\alpha)}
|\theta^{i,j}_{\eps}-c_k|^{p/(p-\alpha)}
m(Y^i_{\eps}\times{T^j_{0,\eps}})\right) \\
\quad\leq{C}\left(\sum_{k=0}^m
\left(m(E^{i,j,k}_{\eps})+\int_{E^{i,j,k}_{\eps}}|\varphi|^p\,dxdt
+|c_k|^p m(E^{i,j,k}_{\eps})\right)\right)^{\frac{p-1-\alpha}{p-\alpha}}
\|\varphi-\Psi\|^{p/(p-\alpha)}_{\cal U} \\
\quad{}+{C}\sum_{k=0}^m\left(
\left(m({\tilde {E}}^{i,j,k}_{\eps})+
\int_{{\tilde {E}}^{i,j,k}_{\eps}}|\varphi|^p\,dxdt
+|c_k|^p m({\tilde {E}}^{i,j,k}_{\eps})\right)^{\frac{p-1-\alpha}{p-\alpha}}
\|\varphi-c_k\|^{p/(p-\alpha)}_{\cal U}\right).
\end{gather*}
Hence
\begin{gather}
\|p_{\eps}(\cdot,\cdot,M_{\eps}\varphi)
-p_{\eps}(\cdot,\cdot,\Psi)\|^p_{\cal U} \notag\\
\quad\leq C(m(\om\times{]0,{\rm T}[})+\|\varphi\|^p_{\cal U}+
\|\Psi\|^p_{\cal U})^{(p-1-\alpha)/(p-\alpha)}
\|\varphi-\Psi\|^{p/(p-\alpha)}_{\cal U} \notag \\
\quad{}+{C}\sum_{k=0}^m\left(
\left(m({\tilde {E}}^{i,j,k}_{\eps})+
\int_{{\tilde {E}}^{i,j,k}_{\eps}}|\varphi|^p\,dxdt
+|c_k|^pm({\tilde {E}}^{i,j,k}_{\eps})\right)^{\frac{p-1-\alpha}{p-\alpha}}
\|\varphi-c_k\|^{p/(p-\alpha)}_{\cal U}\right).
\end{gather}
Now recall that $m(\p\om_k)=m(\p\delta_k)=0$ for $k\neq{0}$. Thus,
$m({\tilde {E}}^{i,j,k}_{\eps})\to{0}$ as $\eps\to{0}$ for
every $k=0,\,1,\,\ldots,\,m$ and the lemma is proved.
\end{proof}

\section{Proof of the main corrector result}

In this section we give the proof of the main corrector result,
Theorem 3.1, stated in Section 3. Our proof will follow
the lines of the proof of the corrector result for the
corresponding elliptic problem, earlier proved by
Dal Maso and Defranceschi in [4]. We start out by proving
an estimate on $p_{\eps}(\cdot,\cdot,M_{\eps}Du)$
uniformly with respect to $\eps$.
\begin{lemma}
Let $p_{\eps}$ be defined as in {\rm (3.5)}. Then,
\begin{equation}
\|p_{\eps}(\cdot,\cdot,M_{\eps}Du)\|^p_{\cal U}\leq{C},
\end{equation}
where the positive constant $C$ is independent of $\eps$.
\end{lemma}
\begin{proof}
Let us define
\[
\theta^{i,j}_{\eps} = {1\over m(Y^i_{\eps}\times{T}^j_{0,\eps})}
\int_{T^j_{0,\eps}}\int_{Y^i_{\eps}}Du(x,t)\,dxdt.
\]
and
\begin{gather*}
{\tilde {I}}_{\eps}=\{i\in{\bf Z}^{\rm N}:
Y^i_{\eps}\cap\om\neq\phi,\;Y^i_{\eps}\setminus\om\neq\phi\},\\
{\tilde {J}}_{\eps}=\{j\in{\bf Z}:
T^i_{0,\eps}\cap{]0,{\rm T}[}\neq\phi,
\;T^i_{0,\eps}\setminus{]0,{\rm T}[}\neq\phi\}.
\end{gather*}
We apply Lemma 4.1, Lemma 4.2 and the inequality (3.3) to obtain
\begin{gather}
\|p_{\eps}(\cdot,\cdot,M_{\eps}Du)\|^p_{\cal U} \notag \\
\quad=\sum_{j\in{J}_{\eps}}\sum_{i\in{I}_{\eps}}
\int_{T^j_{0,\eps}}\int_{Y^i_{\eps}}
|p_{\eps}(x,t,\theta^{i,j}_{\eps})|^p\,dxdt+
\int_{]0,{\rm T}[\setminus\delta_{\eps}}
\int_{\om\setminus\om_{\eps}}|p_{\eps}(x,t,0)|^p\,dxdt \notag \\
\quad\leq\sum_{j\in{J}_{\eps}}\sum_{i\in{I}_{\eps}}
C(1+|\theta^{i,j}_{\eps}|^p)m(Y^i_{\eps}\times{T}^j_{0,\eps}) 
\notag\\
\quad\quad{}+\left(\sum_{j\in{\tilde {J}}_{\eps}}\sum_{i\in{\tilde{I}}_{\eps}}
\|p_{\eps}(\cdot,\cdot,0)\|^{p+\eta}_{L^{p+\eta}(T^j_{0,\eps};
L^{p+\eta}(Y^i_{\eps};{\bf R}^{\rm N}))}\right)^{p/(p+\eta)} \notag \\
\quad\quad{}\times m((\om\setminus\om_{\eps})
\times(]0,{\rm T}[\setminus\delta_{\eps}))^{\eta/(p+\eta)} \notag \\
\quad\leq Cm(\om\times{]0,{\rm T}[})+C\|M_{\eps}Du\|^p_{\cal U}+
C\left(\sum_{j\in{\tilde {J}}_{\eps}}\sum_{i\in{\tilde{I}}_{\eps}}
m(Y^i_{\eps}{}\times{T}^j_{0,\eps})\right)^{p/(p+\eta)} \notag \\
\quad\quad\times m((\om\setminus\om_{\eps})
\times(]0,{\rm T}[\setminus\delta_{\eps}))^{\eta/(p+\eta)} \notag \\
\quad\leq Cm(\om\times{]0,{\rm T}[})+C\|Du\|^p_{\cal U}+
C\left(\sum_{j\in{{\tilde {J}}}_{\eps}}\sum_{i\in{\tilde I}_{\eps}}
m(Y^i_{\eps}\times{T}^j_{0,\eps})\right)^{p/(p+\eta)} \notag \\
\quad\quad{}\times m((\om\setminus\om_{\eps})
\times(]0,{\rm T}[\setminus\delta_{\eps}))^{\eta/(p+\eta)}. 
\end{gather}
Now $\sum_{j\in{\tilde {J}}_{\eps}}\sum_{i\in{\tilde{I}}_{\eps}}
m(Y^i_{\eps}\times{T}^j_{0,\eps})$ approaches $m(\p\om\times\p{]0,{\rm T}[})$
and $m((\om\setminus\om_{\eps})\times(]0,{\rm T}[\setminus\delta_{\eps}))$
tends to zero
as $\eps\to{0}$. Thus, (5.1) follows by (5.2) and Lemma 5.1 is proved.
\end{proof}

\begin{proof}[Proof of Theorem 3.1]
By the strict monotonicity assumption
it follows that
\begin{gather}
\|p_{\eps}(\cdot,\cdot,M_{\eps}Du)-Du_{\eps}\|_{\cal U} \notag \\
\quad\leq C\left(\int^T_0\int_{\om}
(a_{\eps}(x,t,p_{\eps}(x,t,M_{\eps}Du))-a_{\eps}(x,t,Du_{\eps}),
p_{\eps}(x,t,M_{\eps}Du)-Du_{\eps})\,dxdt\right)^{1/p}. 
\end{gather}

Consequently, Theorem 3.1 is proved if we can prove that
\begin{equation}
\int^T_0\int_{\om}
(a_{\eps}(x,t,p_{\eps}(x,t,M_{\eps}Du))-a_{\eps}(x,t,Du_{\eps}),
p_{\eps}(x,t,M_{\eps}Du)-Du_{\eps})\,dxdt\to{0}
\end{equation}
as $\eps\to{0}$. The proof of (5.4) will be splitted up into four steps.\\\\
{\bf Step 1.} We start by showing that
\begin{equation}
\int^T_0\int_{\om}
(a_{\eps}(x,t,p_{\eps}(x,t,M_{\eps}Du)),
p_{\eps}(x,t,M_{\eps}Du))\,dxdt\to
\int^T_0\int_{\om}(b(Du),Du)\,dxdt.
\end{equation}
Let us write
\begin{gather}
\int^T_0\int_{\om}
(a_{\eps}(x,t,p_{\eps}(x,t,M_{\eps}Du)),
p_{\eps}(x,t,M_{\eps}Du))\,dxdt \notag\\
\qquad=\sum_{j\in{J}_{\eps}}\sum_{i\in{I}_{\eps}}
\int_{T^j_{0,\eps}}\int_{Y^i_{\eps}}
(a({x\over\eps},{t\over\eps^\mu},
p({x\over\eps},{t\over\eps^\mu},\theta^{i,j}_{\eps})),
p({x\over\eps},{t\over\eps^\mu},\theta^{i,j}_{\eps}))\,dxdt \notag \\
\qquad\quad{}+\int_{]0,{\rm T}[\setminus\delta_{\eps}}
\int_{\om\setminus\om_{\eps}}
(a_{\eps}(x,t,p_{\eps}(x,t,0)),p_{\eps}(x,t,0))\,dxdt \notag \\
\qquad=\eps^{\rm N+1}\sum_{j\in{J}_{\eps}}\sum_{i\in{I}_{\eps}}
\int_{T_0}\int_{Y}
(a(y,\tau,
p(y,\tau,\theta^{i,j}_{\eps})),
p(y,\tau,\theta^{i,j}_{\eps}))\,dyd\tau \notag \\
\qquad\quad{}+\int_{]0,{\rm T}[\setminus\delta_{\eps}}
\int_{\om\setminus\om_{\eps}}
(a_{\eps}(x,t,p_{\eps}(x,t,0)),p_{\eps}(x,t,0))\,dxdt \notag \\
\qquad=\sum_{j\in{J}_{\eps}}\sum_{i\in{I}_{\eps}}
\int^T_{0}\int_{\om}
\chi_{Y^i_{\eps}}(y)\chi_{T^j_{0,\eps}}(t)(b(\theta^{i,j}_{\eps}),
\theta^{i,j}_{\eps})\,dyd\tau \notag \\
\qquad\quad{}+\int_{]0,{\rm T}[\setminus\delta_{\eps}}
\int_{\om\setminus\om_{\eps}}
(a_{\eps}(x,t,p_{\eps}(x,t,0)),p_{\eps}(x,t,0))\,dxdt,
\end{gather}
where the last equality follows from (3.8). According to
Remark 1 the map $\varphi\to{b(\varphi)}$ is
continuous from ${\cal U}$ into ${\cal U}'$ and an application of
(3.2), using this fact, yields
\begin{equation}
b(M_{\eps}Du)\to{b(Du)}\;\mbox{ strongly in }\;{\cal U}'.
\end{equation}
and, thus,
\begin{gather}
\sum_{j\in{J}_{\eps}}\sum_{i\in{I}_{\eps}}
\int^T_{0}\int_{\om}
\chi_{Y^i_{\eps}}(y)\chi_{T^j_{0,\eps}}(t)(b(\theta^{i,j}_{\eps}),
\theta^{i,j}_{\eps})\,dyd\tau \notag \\
\qquad=\int^T_0\int_{\om}(b(M_{\eps}Du,M_{\eps}Du)\,dyd\tau\to
\int^T_0\int_{\om}(b(Du,Du)\,dyd\tau. 
\end{gather}
By the uniform continuity assumption we have
\begin{gather*}
|\int_{]0,{\rm T}[\setminus\delta_{\eps}}
\int_{\om\setminus\om_{\eps}}
(a_{\eps}(x,t,p_{\eps}(x,t,0)),p_{\eps}(x,t,0))\,dxdt| \\
\qquad\leq C\int_{]0,{\rm T}[\setminus\delta_{\eps}}
\int_{\om\setminus\om_{\eps}}
(1+|p_{\eps}(x,t,0)|)^p\,dxdt \\*
\qquad\quad{}+|\int_{]0,{\rm T}[\setminus\delta_{\eps}}
\int_{\om\setminus\om_{\eps}}
(a_{\eps}(x,t,0),p_{\eps}(x,t,0))\,dxdt| \\
\qquad\leq Cm((\om\setminus\om_{\eps})\times(]0,{\rm T}[\setminus\delta_{\eps}))+
C\int_{]0,{\rm T}[\setminus\delta_{\eps}}
\int_{\om\setminus\om_{\eps}}
|p_{\eps}(x,t,0)|^p\,dxdt \\* 
\qquad\quad{}+C\left(\sum_{j\in{J}_{\eps}}\sum_{i\in{I}_{\eps}}
m(Y^i_{\eps}\times{T}^j_{0,\eps})\right)^{1/p'}
\left(\int_{]0,{\rm T}[\setminus\delta_{\eps}}
\int_{\om\setminus\om_{\eps}}
|p_{\eps}(x,t,0)|^p\,dxdt\right)^{1/p}. 
\end{gather*}
By arguing as in Lemma 5.1 we conclude that
\[
|\int_{]0,{\rm T}[\setminus\delta_{\eps}}
\int_{\om\setminus\om_{\eps}}
(a_{\eps}(x,t,p_{\eps}(x,t,0)),p_{\eps}(x,t,0))\,dxdt|\to{0}.
\]
Thus, by taking (5.6) and (5.8) into account we have shown (5.5).\\\\
{\bf Step 2.} We proceed by showing that
\begin{equation}
\int^T_0\int_{\om}(a_{\eps}(x,t,p_{\eps}(x,t,M_{\eps}Du)),Du_{\eps})\,dxdt\to
\int^T_0\int_{\om}(b(Du),Du)\,dxdt.
\end{equation}

Let $\rho>{0}$ be arbitrary. For $Du\in{\cal U}$ there exists a simple function
\[
\Psi = \sum_{k=1}^mc_k\chi_{\om_k}\chi_{\delta_k},
\]
which satisfies the assumptions in Lemma 4.4, such that
\begin{equation}
\|Du-\Psi\|_{\cal U}\leq\rho.
\end{equation}

We write
\begin{gather}
\int^T_0\int_{\om}(a_{\eps}(x,t,p_{\eps}(x,t,M_{\eps}Du)),Du_{\eps})\,dxdt \notag \\
\qquad=\int^T_0\int_{\om}(a_{\eps}(x,t,p_{\eps}(x,t,\Psi)),Du_{\eps})\,dxdt \notag \\
\qquad\quad{}+\int^T_0\int_{\om}(a_{\eps}(x,t,p_{\eps}(x,t,M_{\eps}Du))
-(a_{\eps}(x,t,p_{\eps}(x,t,\Psi)),Du_{\eps})\,dxdt. 
\end{gather}
It follows, for the first integral on the right hand side, that
\begin{gather}
\int^T_0\int_{\om}(a_{\eps}(x,t,p_{\eps}(x,t,\Psi)),Du_{\eps})\,dxdt \notag \\
\qquad=\sum_{k=0}^m\int_{\delta_k}\int_{\om_k}(a_{\eps}(x,t,p_{\eps}(x,t,c_k)),Du_{\eps})\,dxdt, 
\end{gather}
where $c_0=0$ and where $\om_0$ and $\delta_0$ are defined as in the previous section.
By Lemma 4.2, the functions $p_{\eps}(\cdot,\cdot,c_k))$ are bounded in
$L^{p+\eta}(0,T;L^{p+\eta}(\om;{\bf R}^{\rm N}))$. By the structure conditions
this implies that $a_{\eps}(\cdot,\cdot,p_{\eps}(\cdot,\cdot,c_k)$ is uniformly
bounded in $L^{s}(0,T;L^{s}(\om;{\bf R}^{\rm N}))$ for some $s>{p'}$.
From Proposition 2.2 it further follows that the sequence $(Du_\eps)$ is
bounded in ${\cal U}$. Therefore there exists a
number $\si> 1$ such that
\[
\|(a_\eps(\cdot,\cdot,p_\eps(\cdot,\cdot,c_k)),Du_\eps)\|_{L^\si(]0,T[\times\om)}\leq C
\]
uniformly with respect to $\eps$. Hence, up to a subsequence,
\[
(a_\eps(\cdot,\cdot,p_\eps(\cdot,\cdot,c_k)),Du_\eps)\to g_k
\mbox{ weakly in }L^\si(]0,T[\times\om),
\]
as $\eps\to 0$. By proposition 2.2 we know that
\[
a_\eps(\cdot,\cdot,p_\eps(\cdot,\cdot,c_k))\to b(c_k)\mbox{ weakly in }{\cal U}'.
\]
This enables us to use the compensated compactness result Theorem 2.1 in [9]
and conclude that
\[
(a_\eps(\cdot,\cdot,p_\eps(\cdot,\cdot,c_k)),Du_\eps)\to(b(c_k),Du)
\]
in the sense of distributions. Consequently $g_k=(b(c_k),Du)$ and
\[
\sum_{k=0}^m\int_{\delta_k}\int_{\om_k}(a_{\eps}(x,t,p_{\eps}(x,t,c_k)),Du_{\eps})\,dxdt
\to\sum_{k=0}^m\int_{\delta_k}\int_{\om_k}(b(c_k),Du)\,dxdt.
\]
By using (5.12) this gives
\begin{equation}
\int^T_0\int_{\om}(a_{\eps}(x,t,p_{\eps}(x,t,\Psi)),Du_{\eps})\,dxdt\to
\int^T_0\int_{\om}(b(\Psi),Du)\,dxdt.
\end{equation}
For the second integral on the right hand side of (5.11) we observe
that the growth condition on $a_\eps$ together with the H\"{o}lder
inequality gives 
\begin{gather}
|\int^T_0\int_{\om}(a_{\eps}(x,t,p_{\eps}(x,t,M_{\eps}Du))
-(a_{\eps}(x,t,p_{\eps}(x,t,\Psi))),Du_{\eps})\,dxdt| \notag \\
\qquad\leq C\int^T_0\int_{\om}(1+|p_\eps(x,t,M_\eps Du)|^p+
|p_\eps(x,t,\Psi)|^p)^{(p-1-\alpha)/p} \notag \\
\qquad\quad{}\times|p_\eps(x,t,M_\eps Du)-p_\eps(x,t,\Psi)|^\alpha|Du_\eps|\,dxdt \notag\\
\qquad\leq C(\int^T_0\int_{\om}(1+|p_\eps(x,t,M_\eps Du)|^p+
|p_\eps(x,t,\Psi)|^p)^{(p-1-\alpha)/p}) \notag \\
\qquad\quad{}\times(\int^T_0\int_{\om}|Du_\eps|^p\,dxdt)^{1/p}
|p_\eps(x,t,M_\eps Du)-p_\eps(x,t,\Psi)|^p\,dxdt)^{\alpha/p}. 
\end{gather}
By the Lemmas 5.1 and 4.1 the sequences $(p_\eps(\cdot,\cdot,M_\eps Du))$
and $(p_\eps(\cdot,\cdot,\Psi))$ are bounded
in ${\cal U}$. Therefore, by using (5.10), the last inequality in (5.14)
and Lemma 4.4 gives
\begin{equation}
\limsup_{\eps\to 0}|\int^T_0\int_{\om}(a_{\eps}(x,t,p_{\eps}(x,t,M_{\eps}Du))
-(a_{\eps}(x,t,p_{\eps}(x,t,\Psi))),Du_{\eps})\,dxdt|\leq C\rho^\gamma.
\end{equation}
By taking Remark 1 into account we obtain
\begin{gather*}
|\int^T_0\int_{\om}(b(Du)-b(\Psi),Du)\,dxdt| \\
\qquad\leq C\int^T_0\int_{\om}(1+|Du|^p+|\psi|^p)^{(p-1-\gamma)/p}
|Du-\Psi|^\gamma|Du|\,dxdt.
\end{gather*}
Again using the H\"{o}lder inequality, and (5.10), yields
\begin{gather}
|\int^T_0\int_{\om}(b(Du)-b(\Psi),Du)\,dxdt| \notag\\
\qquad\leq C(\int^T_0\int_{\om}(1+|Du|^p+|\psi|^p)\,dxdt)^{(p-1-\gamma)/p}
(\int^T_0\int_{\om}|Du|^p\,dxdt)^{1/p}\rho^\gamma\leq C\rho^\gamma.
\end{gather}
Thus (5.9) follows by the arbitrariness of $\rho$ and Step 2 is
accomplished.\\\\
{\bf Step 3} We show that
\begin{equation}
\int^T_0\int_{\om}(a_\eps(x,t,Du_\eps)),p_\eps(x,t,M_\eps Du))\,dxdt\to
\int^T_0\int_{\om}(b(Du),Du)\,dxdt.
\end{equation}
Let us fix $\delta>0$ and let $\Psi$ be defined as in Step 2.
We write
\begin{gather}
\int^T_0\int_{\om}(a_{\eps}(x,t,Du_\eps)),p_{\eps}(x,t,M_\eps Du)\,dxdt \notag \\
\qquad=\sum_{k=0}^m\int_{\delta_k}\int_{\om_k}(a_{\eps}(x,t,Du_{\eps}(x,t,c_k)),p_\eps(x,t,c_k))\,dxdt
\notag \\
\qquad\quad{}+\int^T_0\int_{\om}(a_{\eps}(x,t,Du_{\eps}),p_\eps(x,t,M_\eps Du)-p_\eps(x,t,\Psi)\,dxdt.
\end{gather}
By similar arguments as in Step 2 we conclude that
\begin{equation}
\sum_{k=0}^m\int_{\delta_k}\int_{\om_k}(a_{\eps}(x,t,Du_{\eps}(x,t,c_k)),p_\eps(x,t,c_k))\,dxdt\to
\int^T_0\int_{\om}(b(Du),\Psi)\,dxdt.
\end{equation}
It also follows, by the H\"{o}lder inequality, that
\begin{gather*}
|\int^T_0\int_{\om}(a_{\eps}(x,t,Du_{\eps}),p_\eps(x,t,M_\eps Du)-p_\eps(x,t,\Psi))\,dxdt|
\\
\qquad\leq(\int^T_0\int_{\om}|a_{\eps}(x,t,Du_{\eps})|^{p'}\,dxdt)^{1/p'} \\
\qquad\quad{}\times
(\int^T_0\int_{\om}|p_\eps(x,t,M_\eps Du)-p_\eps(x,t,\Psi)|^p\,dxdt)^{1/p}.
\end{gather*}
Therefore, according to Lemma 4.4,
\begin{equation}
\limsup_{\eps\to 0}
|\int^T_0\int_{\om}(a_{\eps}(x,t,Du_{\eps}),p_\eps(x,t,M_\eps Du)-p_\eps(x,t,\Psi))\,dxdt|
\leq C\rho^{1/p-\alpha}.
\end{equation}
(5.17) now follows by an analogous argumentation as in the final lines of Step 2.\\\\
{\bf Step 4} In order to conclude the proof let us show that
\begin{equation}
\int^T_0\int_{\om}(a_\eps(x,t,Du_\eps)),Du_\eps))\,dxdt\to
\int^T_0\int_{\om}(b(Du),Du)\,dxdt.
\end{equation}  
First we observe that
\[
\int^T_0\int_{\om}(a_\eps(x,t,Du_\eps)),Du_\eps))\,dxdt=-\la u'_\eps,u_\eps\ra+
\la f,u_\eps\ra,
\]
or equivalently
\[
\int^T_0\int_{\om}(a_\eps(x,t,Du_\eps)),Du_\eps))\,dxdt=
-{1\over 2}(\|u_\eps(T)\|^2_{L^2(\om)}-\|u_\eps(0)\|^2_{L^2(\om)})+\la f,u_\eps\ra.
\]
Since ${\cal W}_0$ is continuously embedded in $C(0,T;L^2(\om))$ we can pass to the
limit in the right hand side and, consequently,
\[
\int^T_0\int_{\om}(b(Du),Du)\,dxdt=
-{1\over 2}(\|u(T)\|^2_{L^2(\om)}-\|u(0)\|^2_{L^2(\om)})+\la f,u\ra.
\]
By collecting the results from the Steps 1-4 (5.21) follows and the proof is complete. 
\end{proof}
\begin{remark} 
The results of Theorem 3.1 remain valid even for non-homogeneous
or even more general boundary data. This follows from Theorem 6.1 in [9]. We can also
allow oscillating right hand side and initial data, c.f. Theorem 4.1 and Remark 6.1 in [9].
\end{remark}

\subsection*{Acknowledgements}

This work has been
supported by the Swedish Natural Science Research
Council, the Swedish Research Council for Engineering Sciences.

\label{lastpage}


\begin{thebibliography}{99}

\small

\bibitem{4}  Bensoussan A.,  Lions J.L.\ and  Papanicolaou G., 
Asymptotic Analysis for Periodic Structures, North-Holland, 1978.

\bibitem{7}  Brahim-Otsmane S.,  Francfort G.A.\ and  Murat F., 
Correctors for the Homogenization of the Wave and Heat Equations, 
{\em J. Math. pures et appl.}, 1992, V.71, 197--231.

\bibitem{9}  Dall'Aglio A. and  Murat F., 
A Corrector Result for H-Converging Parabolic
Problems with Time Dependent Coefficients, 
{\em Ann. Scu. Norm. Super. Pisa}, 1997, V.25, N~1, 329--373.

\bibitem{14}  Dal Maso G.\ and  Defranceschi A., 
Correctors for the Homogenization of
Monotone Operators, 
{\em Differential and Integral Equations}, 1990, V.3, 1137--1152.

\bibitem{23a} DiBenedetto E.\ and Friedman A., Regularity of Solutions to
Nonlinear Degenerate Parabolic Systems, 
{\em J. reine angew. Math.}, 1984, V.349, 83--128.

\bibitem{23b} Holmbom A., Homogenization of Parabolic
Equations, an Alternative Approach and Some Corrector Type
Results, {\em Appl. of Math.}, 1997,  V.42, N~5, 321--343.

\bibitem{23c} Meyers N.G.\ and Elcrat A., Some Results on Regularity for
Solutions of Non-Linear Elliptic Systems and Quasi-Regular Functions, 
{\em Duke Math. J.}, 1975, V.42, 121--136.

\bibitem{23d} Pankov A., G-convergence and Homogenization of Nonlinear 
Partial Differential Operators,
Mathematics and its Applications 422, Kluwer Publ., 1997.

\bibitem{23e} Svanstedt N., G-Convergence of Parabolic Operators, 
{\em Nonlinear Analysis TMA},
1999, V.36, N~7, 807--843.

\bibitem{24} Zeidler E., Nonlinear Functional Analysis and its
Applications, Volumes II~A and II~B, Springer-Verlag, 1990.
\end{thebibliography}
\end{document}